\documentclass{article}
\usepackage{amsfonts}
\usepackage{amsmath}

\newtheorem{theorem}{Theorem}

\newtheorem{lemma}[theorem]{Lemma}

\input{tcilatex}
\begin{document}

\author{GOKHAN SOYDAN, MUSA DEMIRCI AND ISMAIL NACI CANGUL}
\title{THE DIOPHANTINE EQUATION $x^{2}+11^{m}=y^{n}$}
\date{}
\maketitle

\begin{abstract}
The object of this paper is to give a new proof of all the solutions of the
Diophantine equation $x^{2}+11^{m}=y^{n},$in positive integers $x,~y~$with
odd $m>1$ and $n\geq 3.$
\end{abstract}



\section{\protect\large Introduction}

The history of the Diophantine equation 
\begin{equation}
x^{2}+C=y^{n},~x\geq 1,~y\geq 1,~n\geq 3
\end{equation}%
goes back to 1850's. In 1850, Lebesque $\cite{Lebesque}$ proved that the
equation (1) has no solutions when $C=1$. The title equation is actually a
special case of the Diophantine equation $ay^{2}+by+c=dx^{n},$ where $a,b,c$
and $d$ are integers, $a\neq 0$, $b^{2}-4ac\neq 0,$ $d\neq 0,$ which has
only a finite number of solution in integers $x$ and $y$ when $n\geq 3,$ see 
$\cite{Landau}.$ Cohn $\cite{Cohn1}$, solved (1) for most values of $C$ such
that $1\leq C\leq 100.$ Mignotte and de Weger, in $\cite{Mignotte}$, found
solutions of $x^{2}+74=y^{5}$ and $x^{2}+86=y^{5}$. In $\cite{Bugeaud2}%
,\smallskip $ Bugeaud, Mignotte and Siksek covered the remaining ones.

Sometimes, several mathematicians considered, some variations of (1). For
example in $\cite{Tengley},$ all solutions of $x^{2}+B^{2}=2y^{n}$ for $B\in
\{3,4,...,501\}$ were given, where $n\geq 3$ and $(x,y)=1.$ In $\cite{FSAbu4}%
,$ dealt with the solutions of $x^{2}+C=2y^{n}$ where $n\geq 3,$ $x,~y\in 
\mathbb{Z}
^{+}$,$~(x,y)=1$ and $C$ is a positive integer. In $\cite{FSAbu5},$ the
complete solutions of the equation $px^{2}+q^{2m}=y^{p}$ for $p,q$ are
primes, $p>3.$

Recently, several authors became interested in the case where $C$ is a prime
power. In $\cite{Arif4},$ the solutions of $x^{2}+2^{k}=y^{n}$ have been
given under some conditions. In $\cite{Le1},$ the author verifies a
conjecture of Cohn given in $\cite{Cohn2},$ saying that $x^{2}+2^{k}=y^{n}$
has no solutions with $k>2$ is even, and gave three possible solutions. For $%
p=3$, Muriefah in$~\cite{Arif1}$, considered the equation $x^{2}+3^{m}=y^{n}$
for odd $m$ completely and for even $m$ partly, Luca, in $\cite{Luca3},$
completed the solutions of this equation. For $p=5$, the solutions of the
equation $x^{2}+5^{m}=y^{n}$ have been given $\cite{Arif2}$ and $\cite%
{FSAbu7},$ for odd and even values of $m.$ In $\cite{Arif3}$, the authors
dealt with the general case $x^{2}+q^{2k+1}=y^{n}$ for $q$ odd prime, $q\neq
7~(\func{mod}~8)~$and $n\geq 5$ odd. The same authors obtained several
results for $x^{2}+q^{2k}=y^{n}$ in $\cite{FSAbu2}.$ Luca and Togbe, in $%
\cite{Luca4},$ dealt with the equation $x^{2}+7^{2k}=y^{n}$.

Finally, in some recent papers, more complicated cases where $C$ is a
product of more than one prime powers have been considered. For example, in $%
\cite{Luca1}$, the case $x^{2}+2^{a}3^{b}=y^{n}$; in $\cite{Luca2},$ the
case $x^{2}+2^{a}5^{b}=y^{n};$ in $\cite{FSAbu3},~$the case $%
x^{2}+5^{a}13^{b}=y^{n}~$in $\cite{Luca5},$ the case $x^{2}+2^{\alpha
}5^{\beta }13^{\gamma }=y^{n}$ has been studied. In $\cite{Pink}$, Pink
studied the case $x^{2}+2^{a}3^{b}5^{c}7^{d}=y^{n}.$ A survey of these and
many others can be found in $\cite{FSAbu1}$.

\bigskip Here we continue this study with the equation%
\begin{equation}
x^{2}+11^{m}=y^{n},~n\geq 3,~m>1~\text{odd.}
\end{equation}

Our main result is the following.

\begin{theorem}
Let\textbf{\ }$m$ be odd. Then the Diophantine equation $%
x^{2}+11^{m}=y^{n},~m>1,$ $n\geq 3$ has only one solution in positive
integers $x,$ $y$ and the unique solution is given by $m=6M+3,$ $%
x=9324.11^{3M},$ $y=443.11^{2M}$ and $n=3$.
\end{theorem}

The proof of the theorem is divided in two main cases: $(11,x)=1$ and $11|x$.

It is sufficent to consider the case where $x$ is a positive integer. To
prove the theorem we need the following

\begin{lemma}
(Nagell$\cite{Nagell}$) The equation $11x^{2}+1=y^{n}$ where $n$ is odd
integer $\geq 3$ has no solution in integers $x$ and $y$ for $y$ odd and $%
\geq 1$.
\end{lemma}

\section{\protect\large Proof of Theorem 1.}

Let first $m=2k+1$, where $m>1.$ For $m=1,$ Cohn gave two possible solutions
as $x=4$ and $x=58$. If $x$ is odd, then $y$ is even and we get $%
x^{2}+11^{2k+1}\equiv 4$ $(\func{mod}8)$ but as $y^{n}\equiv 0$ $(\func{mod}%
8)$, this is not possible and we take $x$ even and $y$ odd.

$\mathbf{Case}$ $\mathbf{I:}$ Let $(11,x)=1$. First let $n$ be odd. Without
losing generality assume that $n=p$ is an odd prime. Then to find the
solutions of the equation $x^{2}+11^{2k+1}=y^{p}$, we have to consider two
possibilities, (Teo.1,$\cite{Cohn1}$):

\begin{equation}
\pm x+11^{k}\sqrt{-11}=(a+b\sqrt{-11})^{p}
\end{equation}

\bigskip where $y=a^{2}+11b^{2}$ and%
\begin{equation}
x+11^{k}\sqrt{-11}=(\frac{a+b\sqrt{-11}}{2})^{p},\text{ \ }a\equiv b\equiv
1~(\func{mod}2)
\end{equation}

\bigskip where $y=\frac{a^{2}+11b^{2}}{4}.$

In (3), as $y$ is odd, only one of $a$ and $b$ is odd. Equating imaginary
parts we get%
\begin{equation*}
11^{k}=\dsum_{r=0}^{\frac{p-1}{2}}\left( 
\begin{array}{c}
p \\ 
2r+1%
\end{array}%
\right) a^{p-2r-1}.(-11b^{2})^{r}
\end{equation*}

\bigskip Therefore $b$ must be odd. As the summand is not divisible by $11,$
we must have $b=\pm 11^{k}$. Hence 
\begin{equation}
\pm 1=\dsum_{r=0}^{\frac{p-1}{2}}\left( 
\begin{array}{c}
p \\ 
2r+1%
\end{array}%
\right) a^{p-2r-1}.(-11^{2k+1})^{r}
\end{equation}%
and in this case, both signs are impossible,$\cite{Arif3}.$ Hence (3) has no
solutions.

Now we consider equation (4). By equating imaginary parts, we get%
\begin{equation*}
8.11^{k}=b.(3a^{2}-11b^{2})
\end{equation*}%
If $b=\pm 1$ (4) becomes 
\begin{equation*}
\pm 8.11^{k}=3a^{2}-11.
\end{equation*}%
As the case $k=1$ can be easily eliminated, we may suppose $k>1.$ As the
right hand side is always odd and left hand side is even, this has no
solutions.

If $b=\pm 11^{\lambda },~0<\lambda \leq k$, then in this case (4) becomes%
\begin{equation*}
\pm 8.11^{k}=\pm 11^{\lambda }(3a^{2}-11.11^{2\lambda })
\end{equation*}%
\begin{equation*}
\pm 8.11^{k-\lambda }=3a^{2}-11^{2\lambda +1}
\end{equation*}%
and this is not possible modulo $11$ if $k-\lambda >0$. Therefore $k-\lambda
=0$, and therefore%
\begin{equation*}
\pm 8=3a^{2}-11^{2k+1}
\end{equation*}%
and it can easily be shown that the positive sign is not possible. So we
have 
\begin{equation*}
3a^{2}+8=11^{2k+1}
\end{equation*}%
which has the unique solution $a=\pm 21,~k=1~$and $x=9324~\cite{Bugeaud1}$.

Now if $n$ is even, then it is sufficient to consider the case $n=4$, as in
all other cases, the equation can be reduced to case $n=4$ or $n=p$ is an
odd prime. Hence%
\begin{equation*}
(y^{2}-x)(y^{2}+x)=11^{2k+1}
\end{equation*}%
Since $(x,11)=1,$ we get%
\begin{equation*}
y^{2}-x=1\text{ and }y^{2}+x=11^{2k+1}
\end{equation*}%
and therefore%
\begin{equation*}
2y^{2}=11^{2k+1}+1
\end{equation*}%
which is impossible modulo $11.$

$\mathbf{Case}$ $\mathbf{II:}$ Let $11|x$. Then $11|y.$ Let $x=11^{u}.X$ and 
$y=11^{v}.Y$ where $u,v>0$ and $(11,X)=(11,Y)=1$. Then we have%
\begin{equation*}
11^{2v}X^{2}+11^{2k+1}=11^{nv}.Y^{n}
\end{equation*}%
There are three possible cases:

(i) $2u=$min($2u,2k+1,nv)$. Then by cancelling $11^{2u}$ at each side, we
get 
\begin{equation*}
X^{2}+11^{2(k-u)+1}=11^{nv-2u}.Y^{n}
\end{equation*}%
and considering this equation in modulo $11$, we deduce that $nv-2u=0$, so
that%
\begin{equation*}
X^{2}+11^{2(k-u)+1}=Y^{n}
\end{equation*}%
with $(11,X)=1$. If $k-u=0$, by $\cite{Cohn1},$ there are only two possible
solutions of this equation, $(4,3)$ and $(58,15)$. If $k-u>0$, then as
proved before, the only solution is possible when $k-u=1$ and $n=3$.
Therefore $nv=3v=2u$ which implies that $3|u.$ Let $u=3M$. Then $k=3M+1$ and 
$m=6M+3$. So this equation has a solution only when $m=6M+3$ and this
solution is $X=9324$ and $y=443.$ Hence the solution of our equation is $%
x=9324.11^{u}=9324.11^{3M}$ and $y=443.11^{v}=443.11^{2M}$.

(ii) $2k+1=$min($2u,2k+1,nv)$, then%
\begin{equation*}
11^{2u-2k-1}X^{2}+1=11^{nv-2k-1}.Y^{n}
\end{equation*}%
and considering this equation modulo $11,$ we get $nv-2k-1=0$, so $n$ is odd
and $11(11^{u-k-1}X)^{2}+1=Y^{n},$ by Lemma 2 this equation has no solution.

(iii) $nv=$min($2u,2k+1,nv)$. Then $11^{2u-nv}X^{2}+11^{2k+1-nv}=Y^{n}$ and
this is possible modulo $11$ only if $2u-nv=0$ or $2k+1-nv=0$ and both of
these cases have already been discussed. This concludes the proof.

\begin{tabular}{l}
G\"{o}khan Soydan \\ 
Isiklar Air Force High School \\ 
16039 Bursa, TURKEY \\ 
gsoydan@uludag.edu.tr%
\end{tabular}

\bigskip

\begin{tabular}{l}
Musa Demirci, Ismail Naci Cang\"{u}l \\ 
Department of Mathematics \\ 
Uluda\u{g} University \\ 
16059 Bursa, TURKEY \\ 
mdemirci@uludag.edu.tr, cangul@uludag.edu.tr%
\end{tabular}

\end{document}